\documentclass[journal,transmag]{IEEEtran}
\usepackage{cite}
\usepackage{amsmath,amssymb,amsfonts}
\usepackage{mathtools}
\usepackage{arydshln}
\usepackage{graphicx}
\usepackage{textcomp}
\usepackage{xcolor}
\usepackage{algorithm}
\usepackage{algpseudocode}
\usepackage{physics}

\ifCLASSINFOpdf
\else
\fi
\begin{document}
%
% paper title
% Titles are generally capitalized except for words such as a, an, and, as,
% at, but, by, for, in, nor, of, on, or, the, to and up, which are usually
% not capitalized unless they are the first or last word of the title.
% Linebreaks \\ can be used within to get better formatting as desired.
% Do not put math or special symbols in the title.
\title{Proxy-surface-based fast direct solver for TE-mode scattering problems on distributed memory systems}

% author names and affiliations
% transmag papers use the long conference author name format.

\author{\IEEEauthorblockN{Yasuhiro Matsumoto\IEEEauthorrefmark{1} and
Rio Yokota\IEEEauthorrefmark{2}}%,~\IEEEmembership{Fellow,~IEEE}}
\IEEEauthorblockA{\IEEEauthorrefmark{1}Center for Information Infrastructure,
Institute of Science Tokyo, Japan}
\IEEEauthorblockA{\IEEEauthorrefmark{2}Institute of Integrated Research,
Institute of Science Tokyo, Japan}%
\thanks{Manuscript received December 1, 2012; revised August 26, 2015. 
Corresponding author: Y. Matsumoto (email: matsumoto@cii.isct.ac.jp).}}

% The paper headers
\markboth{Journal of \LaTeX\ Class Files,~Vol.~14, No.~8, August~2015}%
%{Shell \MakeLowercase{\textit{et al.}}: Bare Demo of IEEEtran.cls for IEEE Transactions on Magnetics Journals}
%\markboth{IEEE TRANSACTIONS ON MAGNETICS,~Vol.~xx, No.~x, August~2015}%
{Matsumoto \MakeLowercase{\textit{et al.}}: Fast direct solver for TE-mode scattering on distributed memory systems}
% The only time the second header will appear is for the odd numbered pages
% after the title page when using the twoside option.
% 
% *** Note that you probably will NOT want to include the author's ***
% *** name in the headers of peer review papers.                   ***
% You can use \ifCLASSOPTIONpeerreview for conditional compilation here if
% you desire.

% If you want to put a publisher's ID mark on the page you can do it like
% this:
%\IEEEpubid{0000--0000/00\$00.00~\copyright~2015 IEEE}
% Remember, if you use this you must call \IEEEpubidadjcol in the second
% column for its text to clear the IEEEpubid mark.

% use for special paper notices
%\IEEEspecialpapernotice{(Invited Paper)}

% for Transactions on Magnetics papers, we must declare the abstract and
% index terms PRIOR to the title within the \IEEEtitleabstractindextext
% IEEEtran command as these need to go into the title area created by
% \maketitle.
% As a general rule, do not put math, special symbols or citations
% in the abstract or keywords.
\IEEEtitleabstractindextext{%
\begin{abstract}
This paper describes an MPI/OpenMP hybrid parallelized fast direct solver for the scattering problem of transverse electric (TE)-mode electromagnetic waves.
Because TE-mode scattering can be reduced to the two-dimensional Helmholtz equation, solvers based on the hierarchically semiseparable (HSS) representation are highly attractive due to their high parallel efficiency.
However, as the HSS representation applies low-rank approximations to all off-diagonal blocks, it exhibits poor compatibility with high-order discretization methods.
We developed a fast direct solver with $O(h^3)$ convergence for Helmholtz transmission problems, whereas conventional HSS solvers typically yield only $O(h)$ convergence (where $h$ represents intervals between the quadrature nodes).
It is based on the weakly singular Burton--Miller boundary integral equation and the Nystr\"om method with a one-point correction.
Furthermore, recognizing that matrix component calculation, rather than matrix factorization, dominates the total computational time of HSS-type boundary integral solvers, we introduced a load-balancing method to maximize parallel efficiency.
Numerical results demonstrate that the direct solver achieves high-accuracy convergence and nearly ideal strong and weak scalabilities.
\end{abstract}

% Note that keywords are not normally used for peerreview papers.
}

% make the title area
\maketitle

% To allow for easy dual compilation without having to reenter the
% abstract/keywords data, the \IEEEtitleabstractindextext text will
% not be used in maketitle, but will appear (i.e., to be "transported")
% here as \IEEEdisplaynontitleabstractindextext when the compsoc 
% or transmag modes are not selected <OR> if conference mode is selected 
% - because all conference papers position the abstract like regular
% papers do.
\IEEEdisplaynontitleabstractindextext
% \IEEEdisplaynontitleabstractindextext has no effect when using
% compsoc or transmag under a non-conference mode.

% For peer review papers, you can put extra information on the cover
% page as needed:
% \ifCLASSOPTIONpeerreview
% \begin{center} \bfseries EDICS Category: 3-BBND \end{center}
% \fi
%
% For peerreview papers, this IEEEtran command inserts a page break and
% creates the second title. It will be ignored for other modes.
\IEEEpeerreviewmaketitle

\section{Introduction}
%Numerical solutions of large-scale electromagnetic scattering problems are extremely important
%for applications in physics and engineering.
%When considering bounded domains, the finite difference method and
%the finite element method are often used as numerical solvers.
%However, in the case of unbounded domains,
%radiation conditions must be dealt with,
%and numerical solvers based on boundary integral equations are promising.
The numerical solution of large-scale electromagnetic scattering problems is important in both physics and engineering applications.
For problems defined over bounded domains, numerical techniques such as the finite difference method and the finite element method are effective.
However, to solve problems in unbounded domains, these methods become challenging due to the necessity of handling radiation conditions.
Consequently, numerical solvers based on boundary integral equations are a promising alternative.

The fast multipole method (FMM) is well known as an efficient acceleration technique of boundary integral solvers for scattering problems \cite{chew2001fast}.
Parallelization techniques for FMM in distributed memory systems have also been well studied, for example, \cite{ying2003new, 2019Abduljabbar, liu2024massive}.
The FMM is typically interpreted as a method for speeding up matrix-vector multiplication, and is therefore often coupled with Krylov subspace iterative solvers like GMRES.
Consequently, when solving problems involving many right-hand sides, the required computational time, in the worst-case scenario, scales as a multiple of the number of right-hand sides.
The inherent drawback of using iterative solvers with FMM can be effectively handled by employing fast direct solvers (FDS) \cite{martinsson2019}.

\begin{figure}[!t]
\centering
%\includegraphics[width=0.49\linewidth]{figs/weak_admissibility.pdf}
%\hfill
\includegraphics[width=0.35\linewidth]{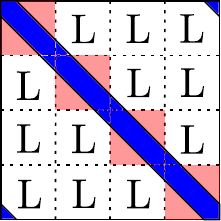}
\caption{Local corrections on the HSS representation.
  The blue band corresponding to the corrected kernel
  spills over into adjacent off-diagonal blocks.
  This figure corresponds to a coefficient matrix stemming from a two-dimensional boundary integral equation when the curvature change of the boundary is not large.
}
\label{fig:weak}
\end{figure}

This study focuses on the transverse electric (TE)-mode electromagnetic scattering in transmission problems.
Because TE-mode scattering reduces to two-dimensional Helmholtz scattering, using solvers based on the hierarchically semiseparable (HSS) representation \cite{chandrasekaran2005fast} is attractive.
FDSs based on the HSS representation include
the Martinsson--Rokhlin solver \cite{MARTINSSON20051},
HSS-ULV factorization \cite{chandrasekaran2006fast}, and
recursive skeletonization \cite{ho2012fast}.
These HSS solvers exhibit high parallel efficiency because all off-diagonal blocks are recursively approximated as low rank, although alternative representations, such as strongly admissible recursive skeletonization \cite{minden2017recursive} and $\mathcal{H}^2$-ULV factorization \cite{ma2024inherently}, are required for three-dimensional problems.
However, the HSS representation is not suitable in high-order discretization because
local corrections or a support of basis functions inevitably spill over into adjacent off-diagonal blocks even when solving two-dimensional problems like that shown in Figure \ref{fig:weak},
leading to a loss of accuracy of low-rank approximation.
Therefore, a major limitation for HSS solvers remains: we are forced to compromise on the discretization order to retain high parallel efficiency.

This study developed an MPI/OpenMP parallelized scalable HSS-type FDS with $O(h^3)$ convergence for TE-mode transmission problems, whereas conventional HSS solvers typically yield only $O(h)$ convergence (where $h$ denotes the quadrature node spacing). This work makes two main contributions: a boundary integral formulation tailored to the HSS representation and an efficient load-balancing technique.

The developed HSS solver is based on the weakly singular Burton--Miller (BM) boundary integral equation, which is expected to avoid quasi-resonances \cite{hiptmair2022spurious}.
Due to this mild singularity, the solver achieves $O(h^3)$ convergence via the Nystr\"om method with a one-point corrected quadrature, in which only the diagonal terms are locally corrected.
Although an efficient parallelization technique for the HSS-ULV factorization has been described \cite{rouet2016distributed},
we observed that the dominant computational bottleneck lies in calculating matrix entries rather than the factorization itself.
These entries are computed during low-rank approximations,
for which we employ the proxy-surface method \cite{MARTINSSON20051}.
While the proxy-surface method has been successfully extended to the Helmholtz \cite{martinsson2007fast}, Maxwell's \cite{rong2019fast}, and Navier--Cauchy \cite{MATSUMOTO2025106148} equations, its parallel load balancing in an HSS solver requires further optimization.
Therefore, we use an appropriate block formulation for the linear system to balance the computational load across processes and threads.
Numerical results demonstrate that our HSS direct solver achieves high-order convergence and nearly ideal strong and weak scalabilities. 
\section{Boundary integral equations for TE-mode transmission scattering problems}
Let us consider the Cartesian coordinate system $(x_1, x_2, x_3)$.
There is a bounded scatterer made of a penetrable object with a smooth surface.
This scatterer is shaped as a general cylinder that extends infinitely along the $x_3$ direction.
A time harmonic electromagnetic wave propagates along the $x_3$ axis.
When the electric field component is entirely perpendicular to the direction of wave propagation,
this specific propagation is called the TE-mode.
Under these assumptions, the scattering problem is reduced to the transmission problem of the Helmholtz equation in the $x_1$--$x_2$ plane.
Let $\Omega_1 \subset \mathbb{R}^2$ be the part of the scatterer cut off at this $x_1$--$x_2$ plane.
We define $\Omega_0 = \mathbb{R}^2 \setminus \Omega_1$ as the exterior connected domain.
For a given function $v$,
$S_{k}v$ and $D_{k}v$ are the single- and double-layer potentials with wavenumber $k$,
and $D_{k}^* v$ and $N_{k} v$ are their normal derivatives, respectively.
To solve the TE-mode scattering, the following boundary integral equation (BIE) can be used:
\begin{multline}
    \mqty(
    -\qty(\varepsilon_{0} S_{k_{0}} + \varepsilon_{1} S_{k_1}) & D_{k_0} + D_{k_1} \\
    -\qty(D_{ k_0}^{\ast} + D_{k_1}^{\ast}) & \frac{1}{\varepsilon_{0}} N_{k_0} + \frac{1}{\varepsilon_{1}} N_{k_1}
    )
    \mqty(
    q(x) \\
    u(x)
    )
    \\ =
    \mqty(
    -u^{\mathrm{in}}(x) \\
    -\frac{1}{\varepsilon_{0}} \pdv{u^{\mathrm{in}}}{n}\qty(x)
    ), \quad x \in \Gamma, \label{eq:pmchwt}
\end{multline}
where $\Gamma = \partial \Omega_1$ is the boundary of a general cylindrical domain in the $x_1$--$x_2$ plane, $n(z)$ is the outward unit normal vector at $z \in \Gamma$,
$u$ is the limit to $\Gamma$ of the magnetic field in the $x_3$ direction,
$u^{\mathrm{in}}$ is the incident wave,
and $q$ is the scaled normal derivative of $u$.
For $j=0, 1$, the wavenumber $k_j$ in $\Omega_j$ is determined by $k_j = \omega \sqrt{\varepsilon_j \mu_j}$,
where $\omega$ is the angular frequency of the incident wave and $\varepsilon_j$ and $\mu_j$ are the constant permittivity and permeability, respectively.
The above BIE, known as the PMCHWT equation, is the standard formulation for transmission problems \cite{chew2001fast}, and it can be discretized using the one-point zeta-corrected quadrature \cite{wu2023}.

As an alternative formulation for the same transmission problem,
the BM equation is given as
\begin{multline}
  \mqty(
  D_{k_1} + I/2 & -\varepsilon_1 S_{k_1} \\
  D_{k_0} - I/2 + \alpha N_{k_0} & -\varepsilon_0 \qty{S_{k_0} +\alpha \qty(D_{k_0}^* + I/2)}
  )
  \\
  \cdot
  \mqty(
  u(x) \\
  q(x) 
  )
  =
  \mqty(
  0 \\
  -u^{\mathrm{in}}(x) - \alpha \pdv{u^{\mathrm{in}}}{n}\qty(x)
  ), \quad x \in \Gamma, \label{eq:bm}
\end{multline}
where $I$ is the identity operator and $\alpha = i/k_0$ is a constant of the BM method.
While the PMCHWT equation \eqref{eq:pmchwt} involves eight layer potentials,
the BM equation contains only six.
This reduction in the number of layer potentials leads to an efficient FDS \cite{matsumoto2026fast}.
However, a solver based on these two equations achieves only $O(h)$ convergence in the one-point zeta-corrected quadrature, because $N_k v$ possesses a hypersingularity, whereas the other layer potentials $S_k v$, $D_k v$, and $D_k^* v$ are only weakly singular.

We therefore employ the weakly singular BM equation, defined as
\begin{multline}
\left(
\begin{array}{c:c}
D_{k_1} + I/2 & -\varepsilon_1 S_{k_1} \\ \hdashline
\begin{aligned}
  D_{k_0} &- I/2 \\
  &+ \alpha \qty(N_{k_0} - N_{k_1}) 
\end{aligned}
&
\begin{aligned}
  -\epsilon_0 S_{k_0} &- \alpha \qty(\epsilon_0 D_{k_0}^* - \epsilon_1 D_{k_1}^* ) \\
  &- \alpha (\varepsilon_0 + \varepsilon_1)I/2
\end{aligned}
\end{array}
\right)
%  \mqty(
%  \begin{pmatrix*}
%  D_{k_1} + I/2 & -\varepsilon_1 S_{k_1} \\ \hdashline
%  \begin{split}
%    D_{k_0} &- I/2 \\
%    &+ \alpha \qty(N_{k_0} - N_{k_1}) 
%  \end{split}
%    &
%  \begin{split}
%  -\epsilon_0 S_{k_0} &- \alpha \qty(\epsilon_0 D_{k_0}^* - \epsilon_1 D_{k_1}^* ) \\
%  &- \alpha (\varepsilon_0 + \varepsilon_1)I/2
%  \end{split}
%%  )
%  \end{pmatrix*}
  \\
  \cdot
  \mqty(
  u(x) \\
  q(x) 
  )
  =
  \mqty(
  0 \\
  -u^{\mathrm{in}}(x) - \alpha \pdv{u^{\mathrm{in}}}{n}\qty(x)
  ), \quad x \in \Gamma. \label{eq:weakly_bm}
\end{multline}
This formulation is equivalent to the one described in \cite{MISAWA201912-191201}.
Because the difference $(N_{k_0} - N_{k_1})v$ has at most a weakly singular kernel,
this achieves $O(h^3)$ convergence even with the one-point corrected quadrature.
We note that, once the boundary values $u$ and $q$ are obtained,
the solution inside the domain can be calculated through integral representations.

\section{Distributed parallel fast direct solver}
This study developed a load-balanced, MPI/OpenMP hybrid parallel FDS.
This solver is a variant of the Martinsson--Rokhlin solver \cite{MATSUMOTO2025106148} and 
features two key characteristics that improve load balancing:
\begin{itemize}
\item The adoption of coarse-grained MPI processes to minimize the volume of MPI communications. Unlike the FMM, which involves the communication of vectors only, this solver requires the exchange of both vectors and matrices between processes.
\item The introduction of a block formulation to ensure that all OpenMP threads handle the same number of layer potentials. As shown in \eqref{eq:weakly_bm}, the number of layer potentials is non-uniform across the $2 \times 2$ blocks.
\end{itemize}
We note that the direct solver follows the standard procedure except for its parallelization scheme.

\begin{figure}[!t]
  \centering
  \includegraphics[width=1.0\linewidth]{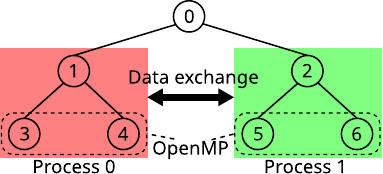}
\caption{Binary tree structure and overview of the parallelization strategy. Each process is assigned specific subtrees of the entire tree. Data required from other processes are exchanged via MPI, while the computational tasks within each subtree are parallelized using OpenMP.}
  \label{fig:tree}
\end{figure}
An overview of the parallelization strategy is shown in Figure \ref{fig:tree}.
First, the boundary $\Gamma$ is recursively partitioned using a perfect binary tree structure.
The entire tree is partitioned into subtrees.
According to this partitioning, the entire index $I_0$ is partitioned into $p$ subsets $I_1, I_2, \ldots, I_p$.
Each MPI process owns a single subtree.
The linear system arising from \eqref{eq:pmchwt}, \eqref{eq:bm}, and \eqref{eq:weakly_bm} can be discretized into the following block forms via partitioning of the index:
\begin{equation}
  \mqty(
  A_{11} & A_{12} & \cdots & A_{1p} \\
  A_{21} & A_{22} & \cdots & A_{2p} \\
  \vdots & \vdots &\ddots &\vdots \\
  A_{p1} & A_{p2} & \cdots & A_{pp}
  )
  \mqty(
  \varphi_1 \\
  \varphi_2 \\
  \vdots \\
  \varphi_p
  )
  =
  \mqty(
  b_1 \\
  b_2 \\
  \vdots \\
  b_p
  ),
  \label{eq:block_linear}
\end{equation}
where $\varphi_j$ and $b_j$ correspond to the discretized solution and right-hand side, respectively.
Each block $A_{ij}$ preserves the original structure of the equation, meaning that it involves the same number of layer potentials.
For example, in the formulation using the weakly singular BM equation \eqref{eq:weakly_bm}, $A_{ij}$ is expressed as
\begin{multline}
  A_{ij} = \\
\left(
\begin{array}{c:c}
D_{k_1}^{ij} + \delta_{ij} I/2 & -\varepsilon_1 S_{k_1}^{ij} \\ \hdashline
\begin{aligned}
  &D_{k_0}^{ij} - \delta_{ij} I/2 \\
  &+ \alpha \qty(N_{k_0}^{ij} - N_{k_1}^{ij})
\end{aligned}
&
\begin{aligned}
  &-\epsilon_0 S_{k_0}^{ij} \\ 
  &- \alpha \qty(\epsilon_0 D_{k_0}^{*ij} - \epsilon_1 D_{k_1}^{*ij} ) \\
  &- \alpha (\varepsilon_0 + \varepsilon_1) \delta_{ij} I/2
\end{aligned}
\end{array}
\right),
%  \mqty(
%D_{k_1}^{ij} + \delta_{ij} I/2 & -\varepsilon_1 S_{k_1}^{ij} \\
%\begin{split}
%  D_{k_0}^{ij} &- \delta_{ij} I/2 \\
%  &+ \alpha \qty(N_{k_0}^{ij} - N_{k_1}^{ij})
%\end{split}
%&
%\begin{split}
%  &-\epsilon_0 S_{k_0}^{ij} \\ 
%  &- \alpha \qty(\epsilon_0 D_{k_0}^{*ij} - \epsilon_1 D_{k_1}^{*ij} ) \\
%  &- \alpha (\varepsilon_0 + \varepsilon_1) \delta_{ij} I/2
%\end{split}
%),
\end{multline}
where $Z_{k}^{ij}$ is the submatrix of a discretized layer potential $Z_{k}$ associated with the index set $I_i$ and $I_j$,
and $\delta_{ij}$ is the Kronecker delta.
Within each MPI process, the task of compressing the block linear system associated with the corresponding subtree is parallelized using OpenMP.
By assigning the calculation of $A_{ij}$ to a single OpenMP thread, the computational load can be evenly distributed. Note that $A_{ij}$ is only calculated around the diagonal at the current level of the tree.
To minimize MPI communication, each process locally stores one-dimensional data, such as global indices and the coordinates of the quadrature points.
Within each assigned subtree, low-rank approximations are computed via the proxy-surface method.
Because this method requires the skeleton indices \cite{MARTINSSON20051} of adjacent cells in the entire tree,
we perform a boundary exchange of these indices at the cells corresponding to the subtree endpoints,
similar to the domain decomposition approach employed in finite difference methods.
%Consequently, the total number of MPI processes employed is restricted to a power of 2.

\section{Numerical experiments}
\begin{figure}[!t]
\centering
\includegraphics[width=0.49\linewidth]{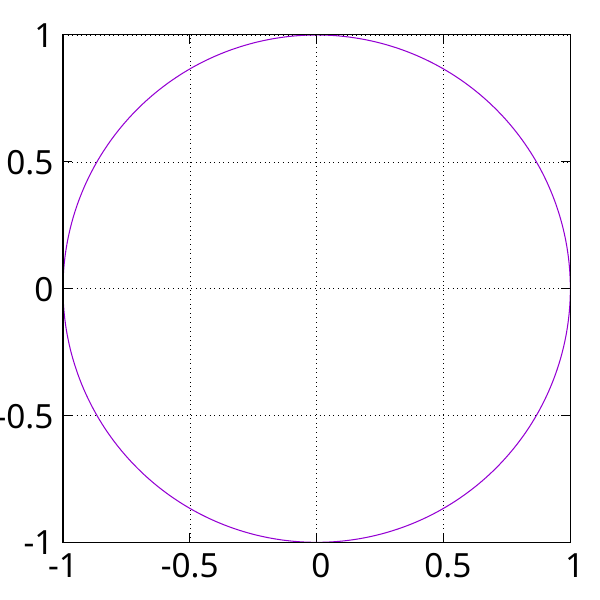}
\hfill
\includegraphics[width=0.49\linewidth]{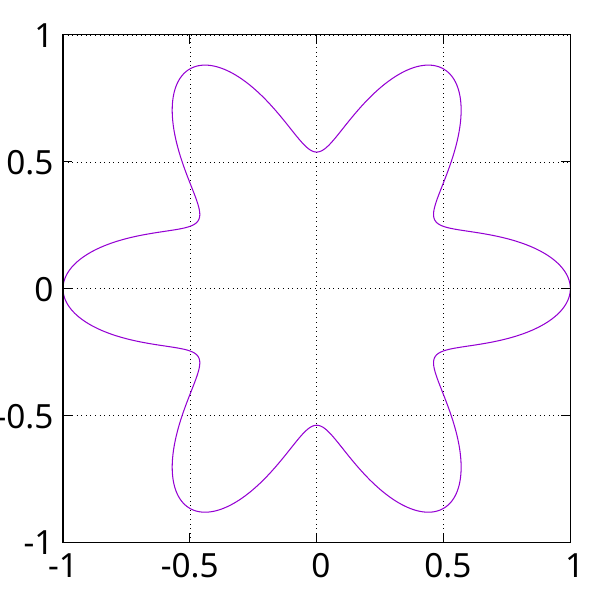}
\caption{Circular- and star-shaped geometries used in the numerical demonstration. Both have a smooth boundary.}
\label{fig:mesh}
\end{figure}

We demonstrate the performance of the direct solver on the TSUBAME4.0 (T4) supercomputer at the Institute of Science Tokyo. As scattering obstacles, we employed the circular and star-shaped geometries illustrated in Fig. \ref{fig:mesh}. For the circular geometry, the analytical solution to the transmission problem can be explicitly calculated. A plane wave propagating in the $x_1$-direction was used as the incident wave. The parameters used in the demonstration are listed in Table \ref{tab:para}. The solver was primarily implemented in C++ using the Eigen (v3) library. In our parallel configuration, a single compute node (comprising 192 CPU cores) of T4 is assigned to each MPI process; thus, the number of nodes is identical to the number of processes.
%While GPUs handle LU factorization for compressed linear equations, CPU-only matrix component computation remains the primary computational bottleneck.

\begin{table}[!t]
%% increase table row spacing, adjust to taste
\renewcommand{\arraystretch}{1.3}
\caption{Parameters used in the numerical demonstration.}
\label{tab:para}
\centering
  \begin{tabular}{ll}
    \hline
    Parameter & Value \\
    \hline
    Shape of scatterer  & Circle, Smooth star \\
    Incident wave       & Plane wave \\
    Material constants & $\varepsilon_{0} = 1$, $\varepsilon_{1} = 5$, $\mu_0 = 1$, $\mu_1 = 1$ \\
    Angular frequency & $\omega = 10$ \\
    \hline
  \end{tabular}
\end{table}

\begin{figure}[!t]
  \centering
  \includegraphics[width=1\linewidth]{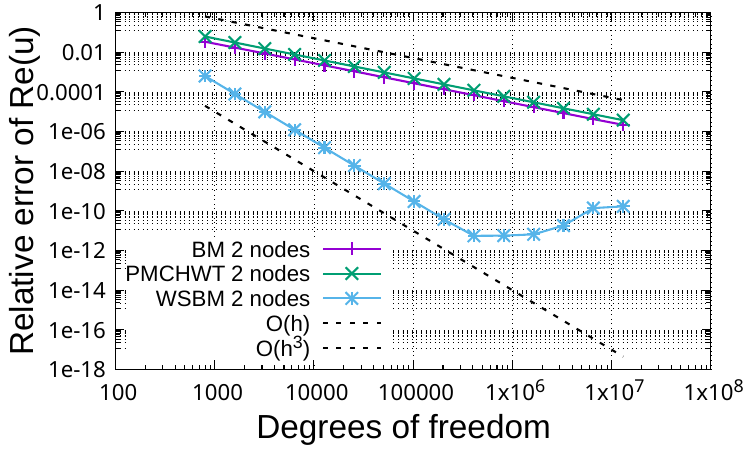}
\caption{Relative 2-norm error of $\Re(u)$ for the solver based on each BIE with a circular geometry.}
\label{fig:error_circle}
\end{figure}

\begin{figure}[!t]
  \centering
  \includegraphics[width=1\linewidth]{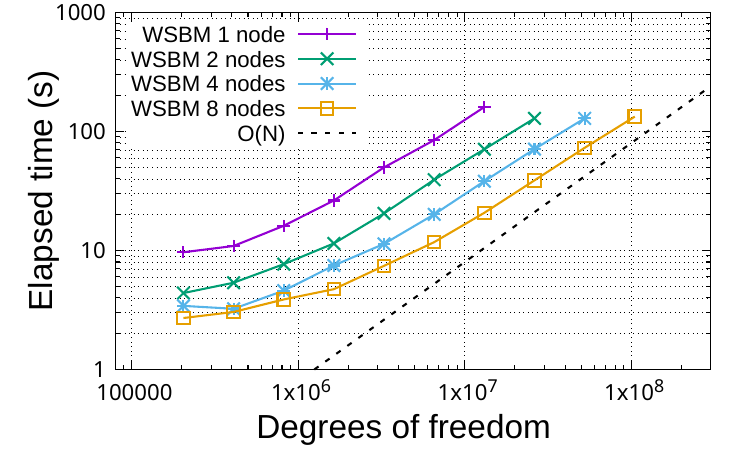}
  \includegraphics[width=1\linewidth]{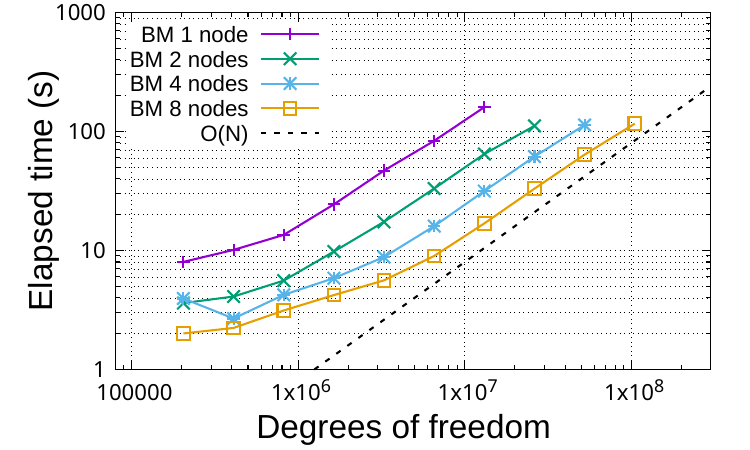}
  \includegraphics[width=1\linewidth]{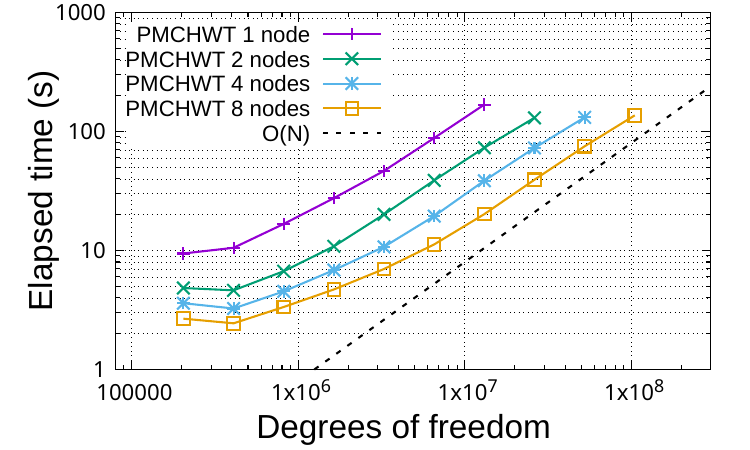}
  \caption{Computational time for the star-shaped geometry using the weakly singular BM equation \eqref{eq:weakly_bm} (top), the standard BM equation \eqref{eq:bm} (middle), and the PMCHWT equation \eqref{eq:pmchwt} (bottom).}   \label{fig:time_star}
\end{figure}

We first verified the implementation for the circular scatterer using two compute nodes.
As shown in Figure \ref{fig:error_circle}, we found that the direct solver based on the weakly singular BM equation \eqref{eq:weakly_bm} achieved $O(h^3)$ convergence, while the solvers based on the other equations achieved only $O(h)$. Here, $h$ denotes the spacing between the quadrature nodes. In this figure and subsequent figures, we refer to the weakly singular BM equation as ``WSBM.''
We next demonstrate the performance of the direct solver for each of the BIEs \eqref{eq:pmchwt}--\eqref{eq:weakly_bm}. Figure \ref{fig:time_star} illustrates the total computational time versus degrees of freedom using 1--8 compute nodes of T4. We observe that all the solvers exhibited nearly ideal weak and strong scalabilities. For example, for approximately $1.31 \times 10^7$ degrees of freedom, the solver based on the weakly singular BM equation achieved a 7.7-fold speedup when using 8 nodes, corresponding to a strong scaling efficiency of 97\%. The computational time for the standard BM equation was approximately 10\%--20\% shorter than that for the weakly singular BM equation. However, from the viewpoint of the convergence rate, the latter remains the most attractive formulation.
We note that the computational time of the solver based on the PMCHWT equation \eqref{eq:pmchwt} was almost identical to that for the weakly singular BM equation \eqref{eq:weakly_bm}, because the number and types of layer potentials are the same.

\section{Conclusion}
This paper described a scalable distributed parallel fast direct solver for the hierarchically semiseparable representation.
The direct solver is based on the weakly singular BM equation with the one-point zeta-corrected quadrature and on the load-balancing technique accomplished by a block formulation of the linear system.
Numerical results demonstrate that our direct solver achieves $O(h^3)$ convergence and nearly ideal strong and weak scalabilities.
As a future direction, it would be interesting to extend the solver to scale for even more practical problems by supporting a larger number of processes, including non-power-of-two configurations, and incorporating GPU acceleration.

\section*{Acknowledgment}
%This work was supported by the JSPS %Japan Society for the Promotion of Science under
% KAKENHI grant number 24K20783. It was carried out using the TSUBAME4.0 supercomputer at the Institute of Science Tokyo, which was provided through the projects ``Joint Usage/Research Center for Interdisciplinary Large-scale Information Infrastructures (JHPCN)'' and ``High Performance Computing Infrastructure (HPCI)'' in Japan (project ID jh260048).

This work was supported by the JSPS KAKENHI grant number 24K20783. It was carried out using the TSUBAME4.0 at the Institute of Science Tokyo, which was provided through the projects JHPCN and HPCI in Japan (project ID jh260048).

% Can use something like this to put references on a page
% by themselves when using endfloat and the captionsoff option.
\ifCLASSOPTIONcaptionsoff
  \newpage
\fi


\begin{thebibliography}{10}
\providecommand{\url}[1]{#1}
\csname url@samestyle\endcsname
\providecommand{\newblock}{\relax}
\providecommand{\bibinfo}[2]{#2}
\providecommand{\BIBentrySTDinterwordspacing}{\spaceskip=0pt\relax}
\providecommand{\BIBentryALTinterwordstretchfactor}{4}
\providecommand{\BIBentryALTinterwordspacing}{\spaceskip=\fontdimen2\font plus
\BIBentryALTinterwordstretchfactor\fontdimen3\font minus
  \fontdimen4\font\relax}
\providecommand{\BIBforeignlanguage}[2]{{%
\expandafter\ifx\csname l@#1\endcsname\relax
\typeout{** WARNING: IEEEtran.bst: No hyphenation pattern has been}%
\typeout{** loaded for the language `#1'. Using the pattern for}%
\typeout{** the default language instead.}%
\else
\language=\csname l@#1\endcsname
\fi
#2}}
\providecommand{\BIBdecl}{\relax}
\BIBdecl

\bibitem{chew2001fast}
W.~C. Chew, E.~Michielssen, J.~Song, and J.-M. Jin, \emph{Fast and efficient
  algorithms in computational electromagnetics}.\hskip 1em plus 0.5em minus
  0.4em\relax USA: Artech House, Inc., 2001.

\bibitem{ying2003new}
L.~Ying, G.~W.  Biros, D.~Zorin, and H.~Langston, ``A new parallel
  kernel-independent fast multipole method,'' in \emph{Proceedings of the 2003
  ACM/IEEE Conference on Supercomputing}, ser. SC '03.\hskip 1em plus 0.5em
  minus 0.4em\relax New York, NY, USA: Association for Computing Machinery,
  2003, p.~14.

\bibitem{2019Abduljabbar}
M.~Abduljabbar, M.~A. Farhan, N.~Al-Harthi, R.~Chen, R.~Yokota, H.~Bagci, and
  D.~Keyes, ``Extreme scale FMM-accelerated boundary integral equation solver
  for wave scattering,'' \emph{SIAM J. Sci. Comput.}, vol.~41, no.~3, pp.
  C245--C268, 2019.

\bibitem{liu2024massive}
X. D. Liu, W. J. He, M. L. Yang, and X. Q. Sheng, ``Massive parallelization of
  multilevel fast multipole algorithm for 3-D electromagnetic scattering
  problems on SW26010 many-core cluster'' \emph{J.
  Supercomput.}, vol.~80, no.~7, pp. 8702--8718, 2024.

\bibitem{martinsson2019}
P. G. Martinsson, \emph{Fast direct solvers for elliptic PDEs}.\hskip 1em plus
  0.5em minus 0.4em\relax Philadelphia, PA: Society for Industrial and Applied
  Mathematics, 2019.

\bibitem{chandrasekaran2005fast}
S.~Chandrasekaran, M.~Gu, and W.~Lyons, ``A fast adaptive solver for
  hierarchically semiseparable representations,'' \emph{Calcolo}, vol.~42,
  no.~3, pp. 171--185, 2005.

\bibitem{MARTINSSON20051}
P. G. Martinsson and V.~Rokhlin, ``A fast direct solver for boundary integral
  equations in two dimensions,'' \emph{J. Comput. Phys.}, vol. 205, no.~1, pp.
  1--23, 2005.

\bibitem{chandrasekaran2006fast}
S.~Chandrasekaran, M.~Gu, and T.~Pals, ``A fast {ULV} decomposition solver for
  hierarchically semiseparable representations,'' \emph{SIAM J. Matrix Anal.
  Appl.}, vol.~28, no.~3, pp. 603--622, 2006.

\bibitem{ho2012fast}
K.~L. Ho and L.~Greengard, ``A fast direct solver for structured linear systems
  by recursive skeletonization,'' \emph{SIAM J. Sci. Comput.}, vol.~34, no.~5,
  pp. A2507--A2532, 2012.

\bibitem{minden2017recursive}
V.~Minden, K.~L. Ho, A.~Damle, and L.~Ying, ``A recursive skeletonization
  factorization based on strong admissibility,'' \emph{Multiscale Model.
  Simul.}, vol.~15, no.~2, pp. 768--796, 2017.

\bibitem{ma2024inherently}
Q.~Ma and R.~Yokota, ``An inherently parallel $\mathcal{H}^2$-{ULV}
  factorization for solving dense linear systems on GPUs,'' \emph{Int. J. High
  Perform. Comput. Appl.}, vol.~38, no.~4, pp. 314--336, 2024.

\bibitem{hiptmair2022spurious}
R.~Hiptmair, A.~Moiola, and E.~A. Spence, ``Spurious quasi-resonances in
  boundary integral equations for the Helmholtz transmission problem,''
  \emph{SIAM J. Appl. Math.}, vol.~82, no.~4, pp. 1446--1469, 2022.

\bibitem{rouet2016distributed}
F. H. Rouet, X.~S. Li, P.~Ghysels, and A.~Napov, ``A distributed-memory package
  for dense hierarchically semi-separable matrix computations using
  randomization,'' \emph{ACM Trans. Math. Software},
  vol.~42, no.~4, pp. 1--35, 2016.

\bibitem{martinsson2007fast}
P. G. Martinsson and V.~Rokhlin, ``A fast direct solver for scattering problems
  involving elongated structures,'' \emph{J. Comput. Phys.}, vol. 221, no.~1,
  pp. 288--302, 2007.

\bibitem{rong2019fast}
Z.~Rong, M.~Jiang, Y.~Chen, L.~Lei, Z.~Nie, and J.~Hu, ``Fast direct surface
  integral equation solution for electromagnetic scattering analysis with
  skeletonization factorization,'' \emph{IEEE Trans. Antennas Propagat.},
  vol.~68, no.~4, pp. 3016--3025, 2019.

\bibitem{MATSUMOTO2025106148}
Y.~Matsumoto and T.~Maruyama, ``Linearly scalable fast direct solver based on
  proxy surface method for two-dimensional elastic wave scattering by cavity,''
  \emph{Eng. Anal. Bound. Elem.}, vol. 173, p. 106148, 2025.

\bibitem{wu2023}
B.~Wu and P. G. Martinsson, ``A unified trapezoidal quadrature method for
  singular and hypersingular boundary integral operators on curved surfaces,''
  \emph{SIAM J. Numer. Anal.}, vol.~61, no.~5, pp. 2182--2208, 2023.

\bibitem{matsumoto2026fast}
Y.~Matsumoto and T.~Maruyama, ``A fast direct solver for two-dimensional
  transmission problems of elastic waves,'' \emph{Eng. Comput.}, vol.~42,
  no.~68, 2026.

\bibitem{MISAWA201912-191201}
R.~Misawa and N.~Nishimura, ``A boundary integral equation suitable for the
  Nystr{\"o}m method having the same complex eigenvalues as the
  {B}urton--{M}iller formulation for the {H}elmholtz equation in 2{D},''
  \emph{Trans. JASCOME}, vol.~19, pp. 61--66, 2019, in Japanese.

\end{thebibliography}
\end{document}